\documentstyle[12pt]{article}
\newcommand{\const}{\mathop{\rm const}\limits}

\newcommand{\supp}{\mathop{\rm supp}\limits}

\newcommand{\diam}{\mathop{\rm diam}\limits}

\newcommand{\Law}{\mathop{\rm Law}\limits}

\newcommand{\Var}{\mathop{\rm Var}\limits}

\newcommand{\Sub}{\mathop{\rm Sub}\limits}

\begin{document}

\begin{center}

{\bf RANDOM  PROCESSES AND  CENTRAL LIMIT THEOREM } \par

\vspace{4mm}

{\bf IN BESOV SPACES }\par

\vspace{4mm}

 $ {\bf E.Ostrovsky^a, \ \ L.Sirota^b } $ \\

\vspace{4mm}

$ ^a $ Corresponding Author. Department of Mathematics and computer science, Bar-Ilan University, 84105, Ramat Gan, Israel.\\

E-mail: eugostrovsky@list.ru\\

\vspace{4mm}

$ ^b $  Department of Mathematics and computer science. Bar-Ilan University,
84105, Ramat Gan, Israel.\\

E-mail: sirota3@bezeqint.net\\

\vspace{5mm}
                    {\sc Abstract.}\\

 \end{center}

 \vspace{4mm}

 \  We  study  sufficient conditions for the belonging of random process  to certain Besov space and for
the Central Limit Theorem (CLT) in these spaces. \par
 \ We investigate also the  non-asymptotic
 tail behavior of normed sums of centered random independent variables (vectors)  with values in these spaces. \par
 \ Main apparatus is the theory of mixed (anisotropic)  Lebesgue-Riesz spaces, in particular so-called permutation
inequality. \par

  \vspace{4mm}

{\it Key words and phrases:} ordinary and generalized Besov spaces, Central Limit Theorem (CLT) in Banach spaces,  mixed (anisotropic)
Lebesgue-Riesz spaces, norms, Orlicz spaces, Monte-Carlo method, Prokhorov's condition,
characteristical functional, Rosenthal constants and inequalities,  metric entropy and entropy integrals, Grand Lebesgue spaces,
exponential  upper tail estimates, triangle and Minkowsky inequalities, permutation inequalities, moments. \par

\vspace{4mm}

{\it 2000 Mathematics Subject Classification. Primary 37B30, 33K55; Secondary 34A34,
65M20, 42B25.} \par

\vspace{4mm}

\section{Notations. Statement of problem.}

\vspace{3mm}

 \  Let $  T  = \{t \} = [0,1] $  be ordinary unit closed numerical segment,   $  f: T \to R  $ be a numerical measurable
function, which may be extended  on the whole  line $  R  $ by zero value:  $ f(t) := 0, \ t \notin T. $ The classical
Lebesgue-Riesz norm  $  |f|_p = |f|_{p,T}  $ of such a function  is defined as follows:

$$
|f|_p = |f|_{p,T} \stackrel{def}{=} \left[  \int_T |f(t)|^p \ dt   \right]^{1/p}, \ p \ge 1; \eqno(1.1)
$$

$$
L_p = L_{p,T} \stackrel{def}{=} \{ f, \ f: T \to R, \ |f|_{p,T} < \infty \}.  \eqno(1.1a)
$$
\ Denote

$$
 S_h[f](t) = S[f](t,h) := f(t + h) - f(t), \ h \in [-1,1], \ t \in [0,1];
$$
and denote also by $ \Delta[f, \delta]_p, \ \delta \in [0,1] $
the $  L_p  $ module of continuity of the function $  f  $ from the space $ L_p:  $

$$
 \Delta[f, \delta]_p \stackrel{def}{=} \sup_{|h| \le \delta} | S_h[f]   |_p. \eqno(1.2)
$$

 It is known that

 $$
 f \in L_{p,T} \ \Rightarrow \ \lim_{\delta \to 0+} \Delta[f, \delta]_p = 0.
 $$

\vspace{3mm}

 We reserve the notation $ \omega(\cdot)  $ for the elements of probability space, which should appear  further. \par

\vspace{3mm}

 Recall that the so-called {\it ordinary}
Besov space $  B^p_{\alpha, s}  $ on the functions defined on considered set $ T = [0,1]  $
may be defined in particular  as a space of all measurable function having a finite norm

$$
|| f ||B^p_{\alpha, s} \stackrel{def}{=} |f|_p +
\left\{ \int_0^1 \left[\delta^{- \alpha} \ \Delta[f,\delta]_p \right]^s \ \frac{d \delta}{\delta}  \right\}^{1/s} =
$$

$$
  |f|_p  + || f ||B^{o,p}_{\alpha, s}, \  p,s = \const, \ge 1, \ \alpha = \const, \eqno(1.3)
$$
where

$$
|| f ||B^{o,p}_{\alpha, s} \stackrel{def}{=}  \left\{ \int_0^1 \left[\delta^{-\alpha}\ \Delta[f,\delta]_p \right]^s \ \frac{d \delta}{\delta}  \right\}^{1/s}.
$$

 \ The detail investigation of these spaces may be found in the classical monographs  \cite{Bennet1}, \cite{Besov1}, \cite{Berg1},
\cite{Peetre1}, \cite{Triebel1}; see also  articles \cite{Feneuil1},  \cite{Gallagher1}. \par

 \ It is offered and considered, in particular, in the  article \cite{Feneuil1} the following generalization of these spaces.
Namely, denote

$$
\Delta_q[f,\delta]_p \stackrel{def}{=} \left[ \int_{-\delta}^{\delta} |S_h[f]_p^q \ dh \right]^{1/q} =
\left[ \int_{-\delta}^{\delta} |S[f](t,h)|_p^q \ dh \right]^{1/q}   , \ q \ge 1, \eqno(1.4)
$$
so that $  \Delta_{\infty}[f,\delta]_p =  \Delta[f,\delta]_p, $ with  evident generalization on the case when $ p = \infty $
or on both the cases $  p = \infty  $ and $ q = \infty.  $\par

 The so-called {\it generalized} Besov space $  B^{p,q}_{\alpha, s}  $ on the functions defined on the considered set $ T = [0,1]  $
may be defined in particular  as a space of all measurable function having a finite norm

$$
|| f ||B^{p,q}_{\alpha, s} \stackrel{def}{=} |f|_p +
\left\{ \int_0^1 \left[\delta^{- \alpha} \ \Delta_q[f,\delta]_p \right]^s \ \frac{d \delta}{\delta}  \right\}^{1/s} =
$$

$$
  |f|_p  + || f ||B^{o,p,q}_{\alpha, s}, \  p,s = \const, \ge 1, \ \alpha = \const, \eqno(1.5)
$$
where the semi-norm $ || f ||B^{o,p,q}_{\alpha, s} $ is defined as follows

$$
|| f ||B^{o,p,q}_{\alpha, s} \stackrel{def}{=}
\left\{ \int_0^1 \left[\delta^{-\alpha}\ \Delta_q[f,\delta]_p \right]^s \ \frac{d \delta}{\delta}  \right\}^{1/s}.
$$

 \ Another approach to the definition of (more generalized) Besov's spaces may be found in articles
 \cite{Gogatishvili1}, \cite{Marcos1}. \par

\vspace{4mm}

 \ {\bf Our goal in this article  is  deducing of sufficient conditions for the belonging of almost all path of random process
to certain Besov space and finding  the sufficient conditions for the Central Limit Theorem (CLT) in these spaces. \par
 \ We intend to investigate also the non-asymptotical tail estimated, for instance, exponential decreasing,
  for the distribution for the norm of normed sums of centered independent random processes
 with values in these spaces.} \par

\vspace{4mm}

 \ The offered here results are formulated in the very simple and {\it natural} terms generated only by the source random process:
metric distance between its values and so on.\par

\vspace{4mm}

  \ The paper is organized as follows. In the next section we consider the case of generalized Besov spaces.
The third section contains the Grand Lebesgue Spaces norm estimation  for Besov norm of random processes.
The case of ordinary Besov spaces is investigated in fourths section. \par
 \ The fifth section is devoted to the Central Limit Theorem in generalized Besov spaces.
The  Central Limit Theorem (CLT) in the classical Besov's spaces is content of the next section. \par
 \ We deduce in the seventh section the non-asymptotical estimates, in particular, exponential ones for the
norm of normed sums of independent random processes with paths in Besov spaces. \par
 \ As ordinary, the last section contains   the concluding remarks, namely, some remarks about possible generalization
 on the multivariate case.\par

\vspace{4mm}

 \ In detail: let $ \xi = \xi(t) = \xi(t,\omega), \ t \in T = [0,1]  $ be numerical bi - measurable
valued random process defined apart from the set $  T  $ on some probability space \\
$ (\Omega = \{\omega\}, B, {\bf P}).  $  Question: under what conditions on the
$  \xi(t) $ almost all its  trajectories belong to certain  Besov space

$$
{\bf P} \left(\xi(\cdot) \in B^{p}_{\alpha, s} \right) = 1 \eqno(1.6)
$$
or more generally satisfies  the Central Limit Theorem (CLT) in these spaces. \par

\vspace{4mm}

 \ This problem for the general separable Banach space $ B  $ instead the Besov's space $ B^{p}_{\alpha, s} $ can be
regarded as a classic, see, e.g.
\cite{Dudley1},  \cite{Ledoux1}, \cite{Ostrovsky1}. The case namely of Besov space is considered in the articles
\cite{Boufoussi1}, \cite{Fageot1}, \cite{Feneuil1}, \cite{Yakovenko1}, \cite{Yakovenko2}. \par

 \ The applications of the CLT in the Banach spaces in the Monte-Carlo method and in statistics is described in the articles
\cite{Frolov1}, \cite{Grigorjeva1}, \cite{Ostrovsky302},  \cite{Ostrovsky303}.\par

 \ We intend to  generalize obtained therein results, for instance, on the case of random processes with exponential
decreasing tails of distributions. \par

\vspace{4mm}

  \ Let us to pay attention that the expression  for $  \Delta_q[f,\delta]_p $ may be written as follows. Denote
 $ T(\delta) = [-\delta, \delta],  $  so that $  T = T(1);   $ then

$$
\Delta_q[f,\delta]_p = | S(\cdot, \cdot) |_{p,T; q,T(\delta)}, \eqno(1.7)
$$
where the right-hand side of an equality (1.7) is a nothing more than a so-called {\it  mixed }, or equally
{\it  anisotropic } Lebesgue-Riesz norm, see \cite{Benedek1},  \cite{Besov1}, chapter 2.\par
 Further,

$$
|| f ||B^{o,p,q}_{\alpha, s}  = |\Delta_q[f, \delta]_p|_{s, T, \mu_{\alpha}}, \eqno(1.8)
$$
where  the measure $  \mu_{\alpha} $ is follow

$$
 \mu_{\alpha}(d \delta) = \delta^{-1 - \alpha s} d \delta, \ \delta \in (0,1). \eqno(1.9)
$$

 \ Therefore, the value $ || f ||B^{o,p,q}_{\alpha, s} $ may be represented also through a three dimensional
mixed norm. \par

\vspace{4mm}

 \ We recall here for readers convenience the definition and used for us simple properties of
 the so-called mixed (anisotropic)  Lebesgue (Lebesgue-Riesz) spaces, which appeared in the famous article of
 Benedek A. and Panzone  R. \cite{Benedek1}.  More detail information about this
spaces with described applications see in the books  of  Besov O.V., Il’in V.P., Nikol’skii S.M.
\cite{Besov1}, chapter 1,2; Leoni G. \cite{Leoni1}, chapter 11; \\ Lieb E., Loss M. \ \cite{Lieb1}, chapter 6. \par

\vspace{3mm}

  Let $ (X_k,A_k,\mu_k), \ k = 1,2,\ldots,l $ be measurable spaces with sigma-finite {\it separable}
non - trivial measures $ \mu_k. $ The separability  denotes that  the metric space
$ A_k  $  relative the distance

$$
\rho_k(D_1, D_2) = \mu_k(D_1 \Delta D_2) = \mu_k(D_1 \setminus D_2) + \mu_k(D_2 \setminus D_1)
$$
is separable.\par
 \ Let also $ p = (p_1, p_2, . . . , p_l) $ be $ l- $ dimensional numerical vector such that \\
$ 1 \le p_j < \infty.$ \par

 \ The anisotropic  (mixed) Lebesgue-Riesz space $ L_{ \vec{p}} $ consists on all the  totally measurable
real valued function  $ f = f(x_1,x_2,\ldots, x_l) = f( \vec{x} ): $

$$
f:  \otimes_{k=1}^l X_k \to R
$$
with finite norm $ |f|_{\vec{p}} =     $

$$
 |f|_{p_1, p_2, \ldots, p_l}  = |f|_{p_1,X_1; \ p_2,X_2; \  \ldots,\ p_l,X_l}=
 |f|_{p_1,X_1, \mu_1; \ p_2,X_2, \mu_2; \ \ldots, \ p_l,X_l, \mu_l}   \stackrel{def}{=}
$$

$$
\left( \int_{X_l} \mu_l(dx_l) \left( \int_{X_{l-1}} \mu_{l-1}(dx_{l-1}) \ldots \left( \int_{X_1}
 |f(\vec{x})|^{p_1} \mu_1(dx_1) \right)^{p_2/p_1 }  \ \right)^{p_3/p_2} \ldots   \right)^{1/p_l}. \eqno(1.10)
$$

 \ In particular, for the r.v. $ \xi  $

$$
  |\xi|_p =  \left[ {\bf E} |\xi|^p \right]^{1/p}, \ p \ge 1.
$$

 \ Note that in general case $ |f|_{p_1,p_2} \ne |f|_{p_2,p_1}, $
but $ |f|_{p,p} = |f|_p. $ \par

 \ Observe also that if $ f(x_1, x_2) = g_1(x_1) \cdot g_2(x_2) $ (condition of factorization), then
$ |f|_{p_1,p_2} = |g_1|_{p_1} \cdot |g_2|_{p_2}, $ (formula of factorization). \par

 \ Note that under conditions of separability of the measures $ \{  \mu_k \} $   these spaces are also  separable  Banach spaces. \par

 \ These spaces appear in the Theory of Approximation, Functional Analysis, theory of Partial Differential Equations,
theory of Random Processes etc. \par

\vspace{3mm}

 \  Let for example $  l = 2; $ we agree to rewrite for clarity the expression for $ |f|_{p_1, p_2}  $ as follows:

$$
|f|_{p_1, p_2} := |  f|_{p_1, X_1; \ p_2, X_2} =  |f|_{p_1, X_1, \mu_1; \ p_2, X_2, \mu_2}.
$$
 \ Analogously,

$$
|f|_{p_1, p_2,p_3} = |  f|_{p_1, X_1; \ p_2, X_2; \ p_3, X_3} = |f|_{p_1, X_1, \mu_1; \ p_2, X_2, \mu_2; \ p_3, X_3, \mu_3}.
$$

 \ Let us give an example. Let $ \eta = \eta(x, \omega )  $ be bi - measurable  random field, $ (X = \{x\}, A, \mu) $ be measurable space,
$ p = \const \in [1,\infty). $ As long as the expectation $ {\bf E } $ is also an integral, we deduce

$$
{\bf E} |\eta(\cdot, \cdot)|^p_{p,X} = {\bf E} \int_X |\eta(x, \cdot)|^p \ \mu(dx)    =
$$

$$
 \int_X  {\bf E}  |\eta(x, \cdot)|^p \ \mu(dx)  = \int_X \mu(dx) \left[ \int_{\Omega} |\eta(x, \omega)|^p \ {\bf P}(d \omega)  \right];
$$

\vspace{3mm}

$$
\left[  {\bf E}|\eta|^m_{p,X} \right]^{1/m} = \left[ \left\{ \int_X |\eta(x)|^p \ \mu(dx) \right\}^{m/p}    \right]^{1/m} =
$$

$$
\left[ \int_{\Omega} \ {\bf P}(d \omega) \left\{ \int_X |\eta(x)|^p \ \mu(dx) \right\}^{m/p}  \right]^{1/m} =
|\eta(\cdot, \cdot)|_{p,X; \ m, \Omega}, \hspace{4mm} p,m \ge 1. \eqno(1.11)
$$

\vspace{3mm}

 \ We will use also the so-called {\it permutation inequality} in the terminology of an article \cite{Adams1};  see also \cite{Besov1},
chapter 1,  p. 24-26. Indeed, let $ (Z, B, \mu) $  be another measurable space  and $ \phi:  (\vec{X},Z) = \vec{X} \otimes Z \to R $ be
measurable function. In what follows
$  \vec{X} = \otimes_k X_k.   $ Let also

$$
  r = \const \ge \overline{p} \stackrel{def}{=} \max_j p_j.
$$

 \  It is true the following inequality (in our notations):

$$
|\phi|_{\vec{p}, \vec{X}; r, Z } \le |\phi|_{r, Z;  \vec{p}, \vec{X}}. \eqno(1.12)
$$

 \ We put in what follows  $  Z = \Omega, \ \mu = {\bf P}.  $

\vspace{4mm}

\section{ The case of generalized Besov spaces.}

\vspace{3mm}

 Let  us introduce some new notations.

$$
V[f](t,z,\delta) := S[f](t, z \cdot \delta), \ t,z,\delta \in T. \eqno(2.1)
$$

  Further, let $ m = \const \ge 1;  $ the Pisier's natural distance  $ \ d_m(t,s), \ t,s \in \ T $ on the set $  T  $
of order $  m, \ m = \const \ge 1 \ $  generated by our random process $  \xi(t) $ is defined as follows

$$
 d_m(t,s) \stackrel{def}{=} |\xi(s) - \xi(t)|_{m, \Omega, {\bf P}} = \left[  {\bf E} |\xi(s) - \xi(t)|^m  \right]^{1/m}, \eqno(2.2)
$$
if there exists (for given value $ m) $ and is finite, see \cite{Pisier1}, \cite{Pisier2}.\par
 So that

$$
|S[\xi] (\cdot, \cdot) |_{m, \Omega} =  |S[\xi] (\cdot, \cdot) |_{m, \Omega,{\bf P}} = d_m(t+h,t). \eqno(2.3)
$$

 Denote also

$$
\sigma_m(h) = \sup_{t \in T} d_m(t+h,t) = \sup_{t \in T} \left[ {\bf E} |\xi(t+h) - \xi(t)|^m \right]^{1/m}. \eqno(2.4)
$$

 Introduce a new measure on at the same set $  T = [0,1] $

$$
\nu(A) = \nu_{\alpha, q,s}(A) = \int_A \delta^{\gamma} \ d \delta, \eqno(2.5)
$$
where

$$
\gamma = \gamma(\alpha, q,s) = \frac{s}{q} - \alpha s - 1. \eqno(2.5a)
$$

\vspace{3mm}

{\bf Theorem 2.1.}  Suppose that for some value $ m \ge \max( p,q,s ) $

$$
| \ V[\xi]  \ |_{m,\Omega, {\bf P}; \ p,T; \ q,T; \ s,T, \nu} < \infty. \eqno(2.6)
$$
 Then

$$
{\bf P} \left( \xi(\cdot) \in B^{o,p,q}_{\alpha, s} \right)  = 1 \eqno(2.7)
$$
and moreover

$$
\left| \ || \xi ||B^{o,p,q}_{\alpha,s} \ \right|_{m,\Omega, {\bf P}} \le  | \ V[\xi]  \ |_{m,\Omega, {\bf P}; \ p,T; \ q,T; \ s,T, \nu}. \eqno(2.7a)
$$

\vspace{3mm}

{\bf Proof.} The expression for $  \Delta_q[\xi,\delta]_p   $ may be rewritten as follows

$$
\Delta_q[\xi,\delta]_p = \left[ \int_{-\delta}^{\delta} \ |S[\xi](t,h)|_{p,T}^q \ dh \right]^{1/q} =
$$

$$
\delta^{1/q} \left[ \int_{-1}^1 |S[\xi](t, \delta z)|_p^q  \ dz \right]^{1/q} =
\delta^{1/q} \left[ \int_{-1}^1 |V[\xi](t, \delta, z)|_p^q  \ dz \right]^{1/q} =
$$

$$
\delta^{1/q} | V(t, \cdot, \cdot) |_{p,T; \ q,T}.
$$

 \ Therefore

$$
|| \xi ||B^{o,p,q}_{\alpha, s}  = |   V [\xi](\cdot, \cdot, \cdot) |_{p,T; \ q, T; \ s,\nu,T}
$$
and correspondingly

$$
\left| \ || \xi ||B^{o,p,q}_{\alpha, s} \ \right|_{m, \Omega,{\bf P}}  =
\left| \   V [\xi](\cdot, \cdot, \cdot) \ \right|_{p,T; \ q, T; \ s,\nu,T}\ |_{m,\Omega, {\bf P}} =
$$

$$
\left| \   V [\xi](\cdot, \cdot, \cdot) \ \right|_{p,T; \ q, T; \ s,\nu,T; \ {m,\Omega, {\bf P}}}.
 \eqno(2.8)
$$
 Since $ m \ge \max(p,q,s),  $ we can use the permutation inequality (1.12):

$$
\left| \ || \xi ||B^{o,p,q}_{\alpha, s} \ \right|_{m, \Omega,{\bf P}} \le
\left| \   V [\xi](\cdot, \cdot, \cdot) \ \right|_{m, \Omega, {\bf P}; \ p,T; \ q, T; \ s,\nu,T} < \infty. \eqno(2.9)
$$
 Thus,

$$
\xi(\cdot) \in B^{o,p,q}_{\alpha, s} \eqno(2.10)
$$
 with probability one, \ Q.E.D. \par

\vspace{4mm}

{\bf Remark 2.1. \ Simplification.} \par

\vspace{3mm}

 \ We intend to convince the reader on the simplicity of conditions of theorem 2.1.
Namely, the expression for the right-hand side of  (2.6) and (2.7a)
may be rewritten as follows.

$$
| \ V[\xi]  \ |_{m,\Omega, {\bf P}; \ p,T; \ q,T; \ s,T, \nu} =
| d_m(t + z \delta, t)|_{p,T; \ q,T; \ s,T,\nu}. \eqno(2.11)
$$

 \ Note that the Pisier's natural distance $  d_m(\cdot, \cdot) $ from (2.2) and therefore in (2.11) may be
relatively easy calculated or estimated in the majority of practical cases. \par

\vspace{3mm}

{\bf Remark 2.2.} \par

\vspace{3mm}

  It is reasonable to choose as a capacity of the value $  m  $ in the theorem 2.1 its minimal value

$$
m := m_0 \stackrel{def}{=}\max(p,q,s).
$$
 Correspondingly

$$
\left( || \xi ||B^{o,p,q}_{\alpha,s} \right)_{m_0,\Omega, {\bf P}} \le
| \ V[\xi]  \ |_{m_0,\Omega, {\bf P}; \ p,T; \ q,T; \ s,T, \nu}. \eqno(2.12)
$$

 But  for exponential non-asymptotical estimations for the tail of distribution of the (random) norm

$$
\tau :=  || \xi ||B^{o,p,q}_{\alpha,s}
$$
 the estimation (2.7a) is very convenient. \par

\vspace{3mm}

{\bf Remark 2.3.} \par

\vspace{3mm}

 If in addition to the conditions of theorem 2.1 for some value $ k \ge p, $ for instance for the value $  k = p $

$$
|\xi(\cdot, \cdot)|_{k, {\bf P}, \Omega; \ p,X} < \infty,  \eqno(2.13)
$$
then

$$
\xi(\cdot) \in B^{p,q}_{\alpha, s} \eqno(2.14)
$$
with probability one and moreover

$$
\left| || \xi ||B^{p,q}_{\alpha,s} \right|_{\min(m,k),\Omega, {\bf P}} \le  | \ V[\xi]  \ |_{m,\Omega, {\bf P}; \ p,T; \ q,T; \ s,T, \nu} +
|\xi(\cdot, \cdot)|_{k, {\bf P}, \Omega; \ p, T}. \eqno(2.14a)
$$

 Indeed, it follows from the condition (2.13) that $  \xi(\cdot) \in L_p(T) $ a.e., see e.g. \cite{Ostrovsky304}.

\vspace{3mm}

 Let us choose for simplicity in (2.14a) $  k = m \ge \max(p,q,s), $ then we obtain for arbitrary such a values $  m  $

$$
\left( || \xi ||B^{p,q}_{\alpha,s} \right)_{m, \Omega, {\bf P}} \le  | \ V[\xi]  \ |_{m,\Omega, {\bf P}; \ p,T; \ q,T; \ s,T, \nu} +
|\xi(\cdot, \cdot)|_{m, {\bf P}, \Omega; \ p,T}. \eqno(2.14b)
$$

\vspace{3mm}

\newpage

{\bf Remark 2.4.} \par

\vspace{3mm}

 An example. Suppose in addition to the the conditions of theorem 2.1

 $$
 d_m(t+ h,t) \le  C_m \ h^{\beta}, \ \beta = \beta(m) = \const \in (0,1].
 $$
 If

$$
\beta + 1/q > \alpha,  \hspace{6mm}  m \ge \max(p,q,s),
$$
then

$$
\xi(\cdot) \in B^{o,p,q}_{\alpha, s}
$$
with probability one and moreover

$$
\left| \ || \xi ||B^{o,p,q}_{\alpha,s} \ \right|_{m,\Omega, {\bf P}} \le  C_m  \ s^{-1/s} \ (\beta - \alpha + 1/q)^{-1/s}.
\eqno(2.15)
$$

\vspace{3mm}

{\bf Remark 2.5.} \par

\vspace{3mm}

 The case when $  p = q = s = \infty  $ correspondent to the so-called H\"older space,
see, e.g. \cite{Ratckauskas1} - \cite{Ratchkauskas4}.\par

\vspace{4mm}

\section{ Grand Lebesgue Spaces norm estimation \\
 for Besov norm of random processes.}

\vspace{3mm}

 \ We recall in the  beginning of this section briefly  the definition and some simple properties
 of the so-called Grand Lebesgue spaces;   more detail
investigation of these spaces see in \cite{Fiorenza3}, \cite{Iwaniec2}, \cite{Kozachenko1}, \cite{Liflyand1}, \cite{Ostrovsky1},
\cite{Ostrovsky2}; see also reference therein.\par

  Recently  appear the so-called Grand Lebesgue Spaces $ GLS = G(\psi) =G\psi =
 G(\psi; B), \  1 < B \le \infty, $ spaces consisting
 on all the random variables (measurable functions)  $ f: \Omega \to R $ with finite norms

$$
   ||f||G(\psi) \stackrel{def}{=} \sup_{p \in (1,B)} \left[ |f|_p /\psi(p) \right]. \eqno(3.1)
$$

  Here $ \psi(\cdot) $ is some continuous positive on the {\it open} interval
$ (1,B) $ function such that

$$
     \inf_{p \in (1,B)} \psi(p) > 0, \ \psi(p) = \infty, \ p > B.
$$
 We will denote
$$
 \supp (\psi) \stackrel{def}{=} (1,B) = \{p: \psi(p) < \infty \}.
$$

 \ The set of all $ \psi $  functions with support $ \supp (\psi)= (1,B) $ will be
denoted by $ \Psi(B) = \Psi(1,B). $ \par
  These spaces are complete Banach spaces and moreover rearrangement invariant, see \cite{Bennet1}, and
are used, for example, in the theory of probability  \cite{Kozachenko1},
  \cite{Ostrovsky1}, \cite{Ostrovsky2}; theory of Partial Differential Equations \cite{Fiorenza3},
  \cite{Iwaniec2};  Functional Analysis \cite{Fiorenza3}, \cite{Iwaniec2},  \cite{Liflyand1},
  \cite{Ostrovsky2}; theory of Fourier series, theory of martingales, mathematical statistics,
  theory of approximation etc.\par

 \ Notice that in the case when $ \psi(\cdot) \in \Psi(\infty)  $ and a function
 $ p \to p \cdot \log \psi(p) $ is convex,  then the correspondent space
$ G\psi(\infty) $ coincides with some {\it exponential} Orlicz space. \par
 \ Conversely, if $ B < \infty, $ then the space $ G\psi(B) $ does  not coincides with
 the classical rearrangement invariant spaces: Orlicz, Lorentz, Marcinkiewicz  etc.\par

  Suppose $   ||f||G\psi \in (0,\infty);  $  then

$$
T_f(u) \stackrel{def}{=} \max ( {\bf P}(f > u),  {\bf P}(f < - u) ) \le \exp \left[ - \tilde{\psi}^*(\ln u)    \right], \ u > e.
$$
 where  $ \tilde{\psi}(p) = p \cdot \ln \psi(p), $ and

$$
 \ g^*(y) = \sup_p (p |y| - g(p))
$$
is Young-Fenchel, or Legendre transform of the function $  g = g(p). $ \par

  \ The last relations implies that in the case when $   ||f||G\psi \in (0,\infty)  $ the function $  f(\cdot) $
obeys the exponential decreasing tail of distribution or equally belongs to some exponential Orlicz's space. \par
 In detail: this $  G\psi  $ space coincides with the {\it exponential} Orlicz's space relative the Young function

$$
N(z) = \exp \left( \tilde{\psi}^*(\ln |z|)  \right) - 1, \  |z| > e,
$$
see \cite{Kozachenko1}. \par

\vspace{3mm}

 Denote for instance

$$
\psi_l(p) = p^{1/l}, \ l = \const > 0.
$$

 \ The case $  \psi(p) = \psi_2(p) = \sqrt{p} = p^{1/2} $  correspondent to the so-called {\it subgaussian space}

$$
\Sub(\Omega)  \stackrel{def}{=} G\psi_2.
$$

\ It is known that there exists  absolute constants  $  C_1, C_2 \in (0, \infty) $ such that  any
 {\it centered} (mean zero) r.v.  $ \eta $ belongs to this space $ \Sub(\Omega) $
with  positive finite norm $  ||\eta||G\psi_2 = \sigma^2 \in (0, \infty)  $  iff

$$
\forall \lambda \in R \ \Rightarrow {\bf E} \exp(\lambda \eta) \le \exp(0.5 \ C_1^2 \ \sigma^2 \ \lambda^2) \eqno(3.2)
$$
or equally

$$
T_f(u)  \le \exp (-C_2 u^2), \ u \ge 0. \eqno(3.2a)
$$

 In the more general case, i.e. when $ ||f||G\psi_l = 1,  $ then

$$
T_f(u) \le \exp \left (- C(l) \ u^l   \right), \ u \ge 0. \eqno(3.2b)
$$

 \ Thus, the theory of Grand Lebesgue spaces allows us to obtain in particular the exponential estimates for tails
of distributions for random variables, or equally estimate  the norm of functions in exponential Orlicz spaces. \par

\vspace{3mm}

{\bf Remark 3.1} If we introduce the {\it discontinuous} function

$$
\psi_{(r)}(p) = 1, \ p = r; \hspace{5mm} \psi_{(r)}(p) = \infty, \ p \ne r, \ p,r \in (1,B)
$$
and define formally  $ C/\infty = 0, \ C = \const \in R^1, $ then  the norm
in the space $ G(\psi_r) $ coincides with the $ L_r $ norm:

$$
||f||G(\psi_{(r)} = |f|_r.
$$

 \ Therefore, the Grand Lebesgue Spaces are the direct generalization  not only  of the
 exponential Orlicz's spaces, but also the classical Lebesgue-Riesz spaces $ L_r. $ \par

\vspace{3mm}

 Let us return to the formulated above problem. We will apply the mentioned in the second section one-dimensional
degree of freedom $  m \ge \max(p,q,s).  $  Denote

$$
\nu(m) =  | \ V[\xi]  \ |_{m,\Omega, {\bf P}; \ p,T; \ q,T; \ s,T, \nu} +
|\xi(\cdot, \cdot)|_{m, {\bf P}, \Omega; \ p,T}. \eqno(3.3)
$$

\vspace{4mm}

 {\bf  Proposition 3.1. }  Suppose in addition to the conditions of theorem 2.1 that for some value
$   B > \max(p,q,s) \hspace{5mm} \Rightarrow \nu(B) < \infty. $  It follows immediately from the assertion of theorem 2.1 (2.14a) that

$$
\left( || \xi ||B^{p,q}_{\alpha,s} \right)G\nu  \le 1. \eqno(3.4)
$$

 As a consequence:

$$
T_{|| \xi ||B^{p,q}_{\alpha,s}}(u) \le \exp \left[ - \tilde{\nu}^*(\ln u)    \right], \ u > e. \eqno(3.4a)
$$

 \vspace{3mm}

  \ Let us consider some examples.  \par

\vspace{3mm}

 {\bf Example 3.1.} \  Suppose under the conditions (and notations) of the remark (2.4) $  C_m \le m^{1/l}, \ l = \const > 0. $
 It follows from the inequality (2.15)

$$
\left| \left| \ || \xi ||B^{o,p,q}_{\alpha,s} \ \right| \right| G\psi_l \le  \ s^{-1/s} \ (\beta - \alpha + 1/q)^{-1/s},
\ \beta > \alpha - 1/q. \eqno(3.5)
$$

 In the case when  in addition $  \xi(t) $ is centered Gaussian process, for instance, is ordinary Brownian motion:  $ \xi(t) = w(t), $
we conclude $ l = 2   $ and following

$$
\left| \left| \ || \xi ||B^{o,p,q}_{\alpha,s} \ \right| \right| \Sub(\Omega) \le \ s^{-1/s} \ (\beta - \alpha + 1/q)^{-1/s}. \eqno(3.5a)
$$
  \ Evidently, for the Brownian motion $  \beta = 1/2. $ \par

\vspace{3mm}

 {\bf Example 3.2.} \ It may happen that in the inequality (2.15) $  C_m = \const = s^{1/s}  $ (for simplicity) and
for the values $  m \ge m_o := \max(p,q,s)  $

$$
\beta = \beta(m) = \alpha - \frac{1}{q} + m^{-\gamma}, \ \gamma = \const > 0. \eqno(3.6)
$$

 We get substituting into (2.15)

$$
\left| \ || \xi ||B^{o,p,q}_{\alpha,s} \ \right|_{m,\Omega, {\bf P}} \le  \ m^{\gamma/s} \eqno(3.7)
$$
at that  the case of the values $  m \in [1,m_0] $  may be considered by means of Lyapunov's inequality.\par

 The estimation (3.7) may be rewritten as follows.

$$
\left| \left| \hspace{3mm} || \xi ||B^{o,p,q}_{\alpha,s} \hspace{3mm} \right| \right| G\psi_{s/\gamma} \le 1  \eqno(3.8)
$$
and thus

$$
T_{|| \xi ||B^{o,p,q}_{\alpha,s} } (u) \le \exp\left(-C(s/\gamma) \ u^{s/\gamma} \right),
 \ C(s/\gamma) \in (0,\infty), \ u > 0. \eqno(3.8a)
$$

\vspace{3mm}

 {\bf Example 3.3.} \ Assume that for some constants $ C \in (0, \infty), \ \theta \ge 2  $

$$
d_m(t, t+h) \le C \ h^{\theta - m}, \ \theta - 1 \le m < \theta.
$$

 Denote

$$
\tilde{\theta} = \theta - \alpha + \frac{1}{q}
$$
and let still $ \ \tilde{\theta} > 1. $\par

 We observe $ \left| \left|\hspace{3mm} \ || \xi ||B^{o,p,q}_{\alpha,s} \hspace{3mm} \right| \right| G\psi^{(\theta)} < \infty, $  where

$$
\psi^{(\theta)} (m) = ( \tilde{\theta} - m)^{-1/s}, \ 1 \le m <  \tilde{\theta},
$$
therefore

$$
T_{|| \xi ||B^{o,p,q}_{\alpha,s} } (u) \le C_1(s,p,q,\theta,\alpha) \ u^{-\tilde{\theta} } \ (\ln u)^{\tilde{\theta}/s}, \ u > e^2.
\eqno(3.9)
$$
 We used some  estimations from an article \cite{Ostrovsky305}. \par

\vspace{4mm}

\section{ The case of ordinary Besov spaces.}

\vspace{3mm}

 \  This case is more complicated. Note first of all

$$
 \Delta[f, \delta]_p \stackrel{def}{=} \sup_{|h| \le \delta} | S_h[f]   |_{p,T}  =
$$

$$
\sup_{|z| \le 1} |\xi(t + z \delta) - \xi(t)|_{p,T} = \sup_{|z| \le 1} \theta(z),
$$
where

$$
\theta(z) = \theta(z, \delta) = |\xi(t + z \delta) - \xi(t)|_{p,T}, \eqno(4.1)
$$
meaning that $  t \in T = [0,1]. $ \par

 Let as before  $  m \ge \max(p,q,s)  $ and denote

$$
\mu(\delta) = \mu_m(\delta):= \sup_{|z| \le 1} |\theta(z,\delta)|_{m,\Omega,{\bf P}; \ p,T}; \eqno(4.2)
$$
and introduce the following distance (more exactly, semi-distance) on the set $  [-1,1]:  $

$$
\rho(z_1, z_2) =  \rho_m(z_1, z_2) :=
 \sup_{\delta \in (0,1)} \left[ \frac{|\theta(z_1,\delta) - \theta(z_2,\delta)|_{m,\Omega,{\bf P}}}{\mu_m(\delta)} \right].
\eqno(4.3)
$$
 \  The finiteness of this distance  for certain segment $  m  \in [\max(p,q,s), B],  \ B = \const \in (\max(p,q,s), \infty] \  $
will be presumed.\par

 \ We intend to apply the Orlicz's spaces norm tail estimates  for the distribution of maximum (supremum) of random fields, based on the
so-called entropy technique, see e.g.  \cite{Fernique1}, \cite{Dudley1}, \cite{Ledoux1},
  \cite{Kozachenko1}, \cite{Ostrovsky1}, chapter 3, sections  3.4 - 3.6, \cite{Pisier2} etc. \par
 \ Some preliminary notations. Denote by $ H(T,d,\epsilon) = H(d,\epsilon)  $ the metric {\it entropy} of the set $ T = [0,1]    $
at the point $  \epsilon > 0 $ relative certain distance function $ d = d(z_1, z_2), \ z_1, z_2 \in T, \ $  i.e. the natural logarithm
of the minimal number of the closed balls with radii $  \epsilon, \ \epsilon > 0  $ in the distance $  d(\cdot, \cdot)  $ which cover
all the set $   T. $   Put also

$$
N(T,d,\epsilon) = N(d,\epsilon)  = \exp H(T,d,\epsilon) = \exp H(d,\epsilon);
$$
and denote for brevity

$$
N_m(\epsilon) := N(T, \rho_m, \epsilon); \hspace{6mm} H_m (\epsilon) := \ln N(T, \rho_m, \epsilon). \eqno(4.4)
$$

 Introduce following G.Pisier \cite{Pisier1},  \cite{Pisier2}; see also \cite{Ostrovsky1}, chapter 3, section 3.17
the variables

$$
V(m)  \stackrel{def}{=} 9 \int_0^{D_m} N_m^{1/m}(\epsilon) \ d \epsilon, \eqno(4.5)
$$
where
$$
D_m = \diam(T,\rho_m) := \sup_{z_1, z_2 \in T} \rho_m(z_1, z_2)\le 2. \eqno(4.5a)
$$

 \ Define the following $ \psi \ - $ function $  \beta = \beta(m):  $

$$
\beta(m) := V(m) \cdot | \mu_m(\cdot) |_{s, T, \nu},
$$
it is meaning $  \delta \in T, $
and suppose the finiteness of these function  for at last one value $  m > \max(p,q,s):  $

$$
\exists L > \max(p,q,s) \ \Rightarrow \beta(L) < \infty.  \eqno(4.6)
$$

\vspace{4mm}

{\bf Theorem 4.1.} Let the condition (4.6) be satisfied. Then almost all the trajectories of  the random process $ \xi(t) $
belong to the space $ B^{o,p,q}_{\alpha,s} $  and herewith

$$
\left| \left| \hspace{3mm} || \xi ||B^{o,p}_{\alpha,s} \hspace{3mm} \right| \right| G\beta \le 1, \eqno(4.7)
$$
with correspondent tail estimation.\par

\vspace{3mm}

{\bf Proof.} Let $  m \in (\max(p,q,s),L). $  Let us consider the normed random process (fields)

$$
\zeta_m(z,\delta) := \frac{\theta(z,\delta)}{\mu_m(\delta)}. \eqno(4.8)
$$
 \ We observe using again the permutation inequality

$$
| \ \zeta_m(z,\delta) \ |_{m, \Omega, {\bf P}} \le 1,
$$

$$
| \ \zeta_m(z_1,\delta) - \zeta_m(z_2,\delta_2) \ |_{m, \Omega, {\bf P}} \le \rho_m(z_1, z_2).
$$
 Since the so-called entropy integral  (4.5) convergent,  one can apply the main result of the article
\cite{Pisier2}; see also \cite{Ostrovsky1}, chapter 3, section 3.17:

$$
\sup_z | \zeta_m(z,\delta) |_{m, \Omega, {\bf P}}  \le V(m),
$$
or equally

$$
\sup_z | \theta(z,\delta) |_{m, \Omega, {\bf P}}  \le V(m) \cdot \mu_m(\delta),
$$
which is equivalent to the estimate(4.7). \par

\vspace{4mm}

{\bf Example 4.1.} Let in (4.7) $  m $ be a fixed number $  m = r \in (\max(p,q,s),  L). $ As long as the ordinary
Lebesgue spaces are the particular case of Grand Lebesgue spaces, we propose from the assertion and conditions of
theorem 4.1 the following $  L_r(\Omega)  $ estimation

$$
 \left| \hspace{3mm} || \xi ||B^{o,p}_{\alpha,s} \hspace{3mm} \right|_{r, \Omega}  \le 1.
$$

\vspace{4mm}

{\bf Example 4.2.}  Let in addition to the conditions of theorem 4.1 in (4.7) $  L = \infty  $ and
$   \beta(m) \le m^{1/l}, \ l = \const > 0.  $ Then

$$
\left| \left| \hspace{3mm} || \xi ||B^{o,p}_{\alpha,s} \hspace{3mm} \right| \right| G\psi_l \le 1,  \eqno(4.9)
$$
with correspondent exponential tail estimation

$$
T_{|| \xi ||B^{o,p}_{\alpha,s} }(u) \le e^{- C(l) u^l }, \ u \ge 1.
$$

\vspace{4mm}

 \ We intend now to offer some important generalization of theorem 4.1. Namely, let $ \lambda_m = \lambda_m(\delta), \ \delta \in [0,1] $
be some family of non-negative continuous functions  such that

$$
\lim_{\delta \to 0+} \frac{\mu_m(\delta)}{\lambda_m(\delta)} = 0, \ \mu_m(\delta) \le \lambda_m(\delta), \ \lambda_m(0+) =
\lambda_m(0) = 0.  \eqno(4.10)
$$

 \ For instance, $   \lambda_m $ may be the $  L_m(\Omega)  $  deterministic component of a factorable module of continuity for the
r.p. $  \xi(t), $ see in detail a preprint \cite{Ostrovsky406}.\par

 Let $  m \in (\max(p,q,s),L), $ and  let us consider the normalized  random process

$$
\tau_m(z,\delta) := \frac{\theta(z,\delta)}{\lambda_m(\delta)}. \eqno(4.11)
$$

 Introduce an another (bounded) distance on the set $  T:  $

$$
r(z_1, z_2) =  r_m(z_1, z_2) :=
 \sup_{\delta \in (0,1)} \left[ \frac{|\theta(z_1,\delta) - \theta(z_2,\delta)|_{m,\Omega,{\bf P}; \ p,T}}{\lambda_m(\delta)} \right].
\eqno(4.12)
$$

\vspace{3mm}

 \ Introduce also again following G.Pisier \cite{Pisier1},  \cite{Pisier2} the variables

$$
V_r(m) = V_{r, \tau}(m) \stackrel{def}{=} 9 \int_0^{D_{m,r}} N^{1/m}(T,  r, \epsilon) \ d \epsilon, \eqno(4.12)
$$

$$
D_{m,r} = \diam(T,\rho_m) := \sup_{z_1, z_2 \in T} r_m(z_1, z_2)\le 2. \eqno(4.12a)
$$

\vspace{3mm}

 \ Define the following new $ \psi \ - $ function $  \beta_r = \beta_r(m):  $

$$
\beta_r(m) := V_r(m) \cdot | \lambda_m(\cdot) |_{s, T, \nu},
$$
it is meaning $  \delta \in T, $
and suppose the finiteness of these function  for at last one value $  m > \max(p,q,s):  $

$$
\exists L_r > \max(p,q,s) \ \Rightarrow \beta_r(L_r) < \infty.  \eqno(4.13)
$$

\vspace{3mm}

\ We observe analogously to the proof of theorem 4.1\\

\vspace{3mm}

{\bf Theorem 4.1a.} Let the condition (4.13) be satisfied. Then almost all the trajectories of  the random process $ \xi(t) $
belong to the space $   B^{o,p}_{\alpha, s} $  and wherein

$$
\left| \left| \hspace{3mm} || \xi ||B^{o,p}_{\alpha,s} \hspace{3mm} \right| \right| G\beta_r \le 1, \eqno(4.14)
$$
with correspondent tail estimation.\par

\vspace{3mm}

 {\bf Example 4.3.}  Let $ \xi(t) = w(t), \ 0 \le t \le1 $ be the ordinary Brownian motion (Wiener's process).  It is well
known that one can  take

 $$
 \lambda_m(\delta) = \sqrt{m} \cdot \sqrt{\delta \cdot |\ln \delta|}, \ 0 \le \delta \le 1/e.\eqno(4.15)
 $$
and $  L_r = \infty. $ \par

 We get from the proposition (4.14) of theorem 4.1a  that if $  \alpha < (2s)^{-1}, $ then almost everywhere

$$
 w(\cdot) \in B^{o,p}_{\alpha, s}
$$
and moreover

$$
{\bf P} \left( || \ w(\cdot)   \ ||B^{o,p}_{\alpha, s} > u \right) \le \exp \left( - C(\alpha,p,s) \ u^2   \right), \
u \ge 1. \eqno(4.16)
$$

\vspace{4mm}

\section{ Central limit theorem in generalized Besov spaces.}

\vspace{3mm}

 \ {\bf 1.}  Let  $  (B, ||\cdot||B )  $  be certain separable Banach space builded on the real valued functions defined on our set $  T $
 and  $ \{ \xi_j \} = \{  \xi_j(t) \}, \ t \in T, \ \xi_1(t) = \xi(t), \ j = 1,2,\ldots  $ be a sequence of {\it centered} in the weak sense:
 $ {\bf E} (\xi_i,b) = 0 \ \forall b \in B^* $ or equally $  {\bf E} \xi_j(t) = 0, \ t \in T $ independent identical distributed
(i.; i.d.)  random variables (r.v.) (or equally random vectors, with at the same abbreviation r.v.)  defined on some non-trivial
probability  space $  (\Omega = \{\omega\}, F, {\bf P})   $ with values  in the space  $ B. $  Denote

 $$
 S_n = S_n(t) = S(n) = S(n,t) = n^{-1/2} \sum_{j=1}^n \xi_j(t), \ n = 1,2,\ldots.
 $$

 If we suppose that the r.v. $ \xi  $ has a weak second moment:

 $$
 \forall b \in B^* \ \Rightarrow  (Rb,b) := {\bf E } (\xi,b)^2 < \infty,
 $$
then the characteristical functional (more exactly,  the sequence of  characteristical functionals)

$$
\phi_{S(n)}(b) := {\bf E} e^{i \ (S(n), b) }
$$
of $ S(n) $ converges as $ n \to \infty $ to the  characteristical functional of (weak, in general case)
Gaussian r.v. $ S = S(\infty)  $ with parameters $ (0,R):   $

$$
\lim_{n \to \infty} \phi_{S(n)}(b) = e^{ - 0,5 (Rb,b) }.
$$
  Symbolically:  $ S \sim N(0,R)  $ or $ \Law(S) = N(0,R).  $ The operator $ R = R_{S} $ is called
 the covariation operator, or variance of the r.v. $  S: $

$$
R = \Var(S);
$$
note that $  R = \Var(\xi). $\par

\vspace{4mm}

 \ We recall the classical definition of the CLT in the space $  B. $ \par
(We will investigate in the sequel the case when the space $  B  $ is our Besov's space $  B = B^{o, p,q}_{\alpha,s}.)  $ \par

 \vspace{3mm}

 \ {\bf Definition 5.1.} {\it  We will say as ordinary that the mean zero r.p. $  \xi(t) $ or equally
  the sequence $  \{ \xi_i  \}, \ \xi(t) = \xi_1(t) $  satisfies the
 CLT in the space $  B, $ write:  $ \{\xi_j \} \in CLT = CLT(B) $ or   simple:  $ \xi \in CLT(B),  $  if the
 limiting Gaussian r.v. $  S $  belongs to the  space $  B $ with probability one:  $  {\bf P} (S \in B) = 1 $ and the sequence
 of distributions $\Law(S(n))  $  converges weakly as $ n \to \infty, $   i.e. in the Prokhorov - Skorokhod sense,
  to the distribution of the r.v. $ S = S(\infty):  $  }\par

$$
\lim_{n \to \infty} \Law(S(n)) = \Law(S). \eqno(5.1)
$$

 \ The  equality (5.1) implies that for any  continuous functional $ F: B \to R $

$$
\lim_{n \to \infty} {\bf  P} ( F(S(n)) < x  ) = {\bf P} (F(S) < x) \eqno(5.2)
$$
 for all positive values $ x. $  \par

 \ In particular,

$$
\lim_{n \to \infty} {\bf  P} ( ||S(n)||B < x  ) = {\bf P} (||S||B < x), \ x > 0.
$$

 \vspace{3mm}

{\bf 2.} The problem of describing of necessary (sufficient) conditions  for the infinite - dimensional CLT in Banach space $  B $
has a long history; see, for instance, the  monographs \cite{Araujo1},  \cite{Dudley1},  \cite{Gine2},
\cite{Ledoux1}, \cite{Ostrovsky1} and articles \cite{Gine1},
\cite{Garling1}, \cite{Zinn1}; see also reference therein.\par
 The applications of considered theorem in statistics and method Monte-Carlo  see, e.g. in
\cite{Frolov1}, \cite{Ostrovsky303}, \cite{Ostrovsky407}, \cite{Ostrovsky408}.\par

\vspace{3mm}

 \ {\bf 3.} The cornerstone of  this problem is to establish the {\it weak compactness } of the distributions generated
in the space $  B  $ by the sequence $  \{  S(n) \}: $

$$
\nu_n(D) = {\bf P} ( S(n) \in D),
$$
where $  D  $ is Borelian set in $  B; $ see \cite{Prokhorov1}; \cite{Billingsley1},  \cite{Billingsley2}. \par

\vspace{3mm}

 \ {\bf 4.} We will apply the famous Rosenthal's constants and inequality,  see the classical work of H.P.Rosenthal
 \cite{Rosenthal1}; see also \cite{Ibragimov1}, \cite{Johnson1},  \cite{Ostrovsky502}, \cite{Pinelis1} etc. \par

\vspace{3mm}

 \ Let  $  p = \const \ge 1, \hspace{4mm}  \{ \zeta_k \} $ be a sequence of numerical centered, i.; i.d. r.v.  with finite
 $ p^{th} $ moment $ | \zeta|_p < \infty. $  The following constants,  more precisely, functions on $ p, $ are called
 constants of Rosenthal-Dharmadhikari-Jogdeo-Johnson-Schechtman-Zinn-Latala-Ibragimov-Pinelis-Sharachmedov-Talagrand-Utev...:

$$
K_R(p) \stackrel{def}{=} \sup_{n \ge 1} \sup_{ \{\zeta_k\} } \left[ \frac{|n^{-1/2} \sum_{k=1}^n \zeta_k|_p}{|\zeta_1|_p} \right].
$$

  \ We will  use  the following ultimate up to an error value $ 0.5\cdot 10^{-5} $  estimate  for $ K_R(p), $ see \cite{Ostrovsky502}
and reference therein:

 $$
 K_R(p) \le \frac{C_R \ p}{ e \cdot \log p}, \hspace{5mm}  C_R = \const := 1.77638.
 $$

  \ Note that for the symmetrical distributed r.v. $ \zeta_k $ the constant $  C_R $ may be reduced  up to a value $ 1.53572 $
and that both the boundaries are  exact.\par

\ {\bf 5.}  We retain here all the notations  (and conditions) of the second section, for instance,
the Pisier's notations \cite{Pisier1},  \cite{Pisier2}

$$
V(m)  \stackrel{def}{=} 9 \int_0^{D_m} N_m^{1/m}(\epsilon) \ d \epsilon, \eqno(5.3)
$$

$$
D_m = \diam(T,\rho_m) := \sup_{z_1, z_2 \in T} \rho_m(z_1, z_2)\le 2.
$$

$$
\beta(m) := V(m) \cdot | \mu_m(\cdot) |_{s, T, \nu}, \eqno(5.4)
$$
and add some news:

$$
\tilde{\beta}(m) :=  K_R(m) \cdot V(m) \cdot | \mu_m(\cdot) |_{s, T, \nu}, \eqno(5.4a)
$$
it is meaning $  \delta \in T, $
and suppose the finiteness of these function  for at last one value $  m > m' \stackrel{def}{=} \max(2,p,q,s):  $

$$
\exists L' > \max(2,p,q,s) \ \Rightarrow \tilde{\beta}(L') < \infty.  \eqno(5.5)
$$

\vspace{4mm}

{\bf 6. \ Theorem 5.1.}  Suppose that for some value $ m \ge \max(2, p,q,s ) $

$$
| \ V[\xi]  \ |_{m,\Omega, {\bf P}; \ p,T; \ q,T; \ s,T, \nu} < \infty. \eqno(5.6)
$$
 Then

$$
 \xi(\cdot) \in CLT \left[B^{o,p,q}_{\alpha, s} \right]  \eqno(5.7)
$$
and moreover

$$
\sup_n \left| \ || \ S_n \ ||B^{o,p,q}_{\alpha,s} \ \right|_{m,\Omega, {\bf P}} \le  \tilde{\beta}(m)
\eqno(5.8)
$$
or equally

$$
\sup_n \left| \left| \ || \ S_n \ ||B^{o,p,q}_{\alpha,s} \ \right| \right|G  \tilde{\beta} \le 1. \eqno(5.8a)
$$

\vspace{3mm}

{\bf 7. Proof.} \\

\vspace{3mm}

 \ We  deduce  using Rosenthal's inequality

$$
\overline{d}_m(t,s)  \le K_R(m) \cdot d_m(t,s),
$$
since here $  m \ge 2. $\par

 We can apply the proposition (2.7a) of theorem 2.1

$$
\sup_n \left| \ || S_n(\cdot) ||B^{o,p,q}_{\alpha,s} \ \right|_{m,\Omega, {\bf P}} \le
 K_R(m) \cdot | d_m(t + z \delta, t)|_{p,T; \ q,T; \ s,T,\nu} < \infty.
$$
 The right - hand of the last  estimate  meaning  that $  t \in T, \ z \in T, \delta \in T.   $ \\

\vspace{4mm}

{\bf 8.}    As long as the Banach space $ B^{o,p,q}_{\alpha,s} $  is separable
and the function $  y \to |y|^m, \ m \ge 1  $ satisfies the well known $ \Delta_2  $ condition, there exists a linear compact \\ operator
$  U: \ B^{o,p,q}_{\alpha,s}  \to  B^{o,p,q}_{\alpha,s},  $  which dependent only on the distribution $ \Law(\xi), \ \xi = \xi_1,  $
such that

$$
{\bf P} \left(U^{-1} \xi \in  B^{o,p,q}_{\alpha,s} \right) = 1 \eqno(5.9)
$$
and moreover

$$
{\bf E} || U^{-1} \xi ||^m B^{o,p,q}_{\alpha,s} < \infty, \eqno(5.10)
$$
 \cite{Ostrovsky2}; see also \cite{Buldygin1}, \cite{Ostrovsky603}.\par

{\bf 9.} \ Let us consider the sequence of the  r.v. in the space $  B^{o,p,q}_{\alpha,s}: \ \eta_k(x) = U^{-1} [\xi_k](x); $ it is also a
 sequence of i., i.d. r.v. in the space $   B^{o,p,q}_{\alpha,s}, $  and we can apply the inequality (5.8) taking the value $  m: $

$$
\sup_n {\bf E} || U^{-1}[ S_n]||^m B^{o,p,q}_{\alpha,s} \le K_R^{m} (p) \ {\bf E} | U^{-1}[\xi]|^m B^{o,p,q}_{\alpha,s} = C_m(p) < \infty. \eqno(5.11)
$$
 We get using Tchebychev's  inequality

 $$
\sup_n {\bf P} \left( ||U^{-1}[ S_n]|| B^{o,p,q}_{\alpha,s} > Z   \right) \le C_m(p)/Z^m < \epsilon,\eqno(5.12)
 $$
for sufficiently greatest values $ Z = Z(\epsilon), \ \epsilon \in (0,1). $ \par
  Denote by $  W  = W(Z) $ the set

$$
W = \{ f: f \in  B^{o,p,q}_{\alpha,s}, \  ||U^{-1}[f]|| B^{o,p,q}_{\alpha,s}  \le Z \}. \eqno(5.13)
$$
 Since the operator $ U $ is compact, the set $  W  = W(Z) $ is compact set in the space $  B^{o,p,q}_{\alpha,s}. $ It follows
from the inequality (5.12) that

$$
\sup_n {\bf P} \left( S(n) \notin W(Z) \right) \le \epsilon.
$$
  Thus, the sequence $ \{  S_n \}  $ satisfies the famous Prokhorov's criterion \cite{Prokhorov1} for  weak compactness of
the family of distributions in the separable metric spaces. \par
 This completes the proof of theorem 5.1. \par

\vspace{4mm}

\section{ Central limit theorem in the classical Besov spaces.}

\vspace{3mm}

 \ Denote

$$
S_n(t) = n^{-1/2} \sum_{i=1}^n \xi_i(t), \hspace{5mm} t \in T = [0,1],
$$
where  the random processes $  \xi_i(t) $ are independent copies of the {\it centered (mean zero)} r.p. $ \xi(t)   $
belonging to the space $ B^{o,p}_{\alpha,s}   $ with probability one.  For instance, $  \xi(\cdot)  $ may satisfy
the conditions of theorem 4.1. However, the conditions of offered further  theorem 6.1 "absorb"  ones in theorem 4.1. \par

 \ Define as above
$$
\theta_n(z) = \theta_n(z, \delta) = |S_n(t + z \delta) - S_n(t)|_{p,T}, \eqno(6.0)
$$
meaning that $  t \in T = [0,1]. $ \par

{\it We assume in the sequel  }

 $$
  m \ge \max(2,p,q,s)  \eqno(6.1)
 $$

 and denote

$$
\overline{\mu}(\delta) = \overline{\mu}_m(\delta):= \sup_{|z| \le 1} |\theta_n(z,\delta)|_{m,\Omega,{\bf P}; \ p,T}; \eqno(6.2)
$$
 It follows from the Rosenthal's inequality

$$
 \overline{\mu}_m(\delta)  \le K_R(m) \cdot   \mu_m(\delta). \eqno(6.3)
$$

 \ Introduce also the following distance (more exactly, semi-distance) on the set $  [-1,1]:  $

$$
 \overline{\rho}(z_1, z_2) = \overline{ \rho}_m(z_1, z_2) := \sup_n \
 \sup_{\delta \in (0,1)} \left[ \frac{|\theta_n(z_1,\delta) - \theta_n(z_2,\delta)|_{m,\Omega,{\bf P}; \ p,T}}{\overline{\mu}_m(\delta)} \right].
\eqno(6.4)
$$

 \ Evidently,

$$
\overline{ \rho}_m(z_1, z_2)  \le K_R(m) \cdot  \rho_m(z_1, z_2). \eqno(6.5)
$$

 \  The finiteness of the distance $   \rho_m(z_1, z_2)   $ for certain segment
 $  m  \in [\max(2, p,q,s), B),  \ B = \const \in (\max(2,p,q,s), \infty] \  $ will be presumed.\par

 \ Introduce again following G.Pisier \cite{Pisier1},  \cite{Pisier2}; see also \cite{Ostrovsky1}, chapter 3, section 3.17
the variables

$$
\overline{V}(m)  \stackrel{def}{=} 9 \int_0^{\overline{D}_m}  N(T, \overline{\rho}_m, \epsilon) \ d \epsilon, \eqno(6.6)
$$
where
$$
\overline{D}_m = \diam(T, \overline{\rho}_m) := \sup_{z_1, z_2 \in T} \overline{\rho}_m(z_1, z_2) < \infty. \eqno(6.7)
$$
Define the following $ \psi \ - $ function $ \overline{ \beta} = \overline{\beta}(m):  $

$$
\overline{\beta}(m) := \overline{V(m)} \cdot | \overline{\mu}_m(\cdot) |_{s, T, \nu},
$$
it is meaning $  \delta \in T; $ and suppose the finiteness of these function  for at last one value $  m > \max(2,p,q,s):  $

$$
\exists L > \max(2,p,q,s) \ \Rightarrow \overline{\beta}(L) < \infty.  \eqno(6.8)
$$

\vspace{4mm}

{\bf Theorem 6.1.} Let the condition (6.8) be satisfied. Then the (centered) r.p. $   \xi(t)   $ satisfies the CLT
in the space $ B^{o,p}_{\alpha,s} $ and moreover

$$
\sup_n \left| \left| \hspace{3mm} || S_n(\cdot) ||B^{o,p}_{\alpha,s} \hspace{3mm} \right| \right| G \overline{\beta} \le 1, \eqno(6.9)
$$
with correspondent tail estimation.\par

\vspace{3mm}

{\bf Proof.} Let $  m \in (\max(2,p,q,s),L). $  We use the propositions  (2.7); (2.7a) of theorem 2.1 applied for the
sequence random processes  $  S_n(t);  $ note that all our  estimations are uniformly in $  n.$ This dives us
the estimate (6.9). \par

 The remainder part  of proof theorem 6.1 is completely analogous to one in the theorem 5.1. \par

\vspace{4mm}

{\bf Example 6.1.} Let in (6.1) (and further) $  m $ be a fixed number $  m = r \in (\max(2,p,q,s),  L). $ As long as the ordinary
Lebesgue spaces are the particular case of Grand Lebesgue spaces, we propose from the assertion and conditions of
theorems 4.1 and 6.1 that the r.p. $  \xi(t) $ satisfies the CLT in the space $  B^{o,p}_{\alpha,s}  $
and there holds the following $  L_r(\Omega)  $ estimation

$$
\sup_n \left| \hspace{3mm} || S_n(\cdot) ||B^{o,p}_{\alpha,s} \hspace{3mm} \right|_{r, \Omega}  < \infty. \eqno(6.10)
$$

\vspace{4mm}

{\bf Example 6.2.}  Let in addition to the conditions of theorem 6.1  $  L = \infty  $ and
$  \overline{\beta}(m) \le m^{1/l}, \ l = \const > 0.  $ Then
the r.p. $  \xi(t) $  satisfies again the CLT in the space $  B^{o,p}_{\alpha,s}  $ and moreover

$$
\sup_n \left| \left| \hspace{3mm} || S_n ||B^{o,p}_{\alpha,s} \hspace{3mm} \right| \right| G\psi_{l/(l+1)} \le 1,  \eqno(6.11)
$$
with correspondent exponential tail estimation

$$
\sup_n T_{||S_n(\cdot) ||B^{o,p}_{\alpha,s} }(u) \le e^{- C_2(l) u^{l/( l + 1  ) }}, \ u \ge 1. \eqno(6.11a)
$$

\vspace{4mm}

 \ It is no hard to generalize this example on the case which was considered in (4.10)-(4.11), theorem 4.1a.
 Indeed, let $  \lambda_m = \lambda_m(\delta), \ \delta \in [0,1]  $ be
some family of non-negative continuous functions  such that

$$
\lim_{\delta \to 0+} \frac{\mu_m(\delta)}{\lambda_m(\delta)} = 0, \ \mu_m(\delta) \le \lambda_m(\delta), \ \lambda_m(0+) =
\lambda_m(0) = 0.  \eqno(6.12)
$$

\vspace{3mm}

 \ Define the following new $ \psi \ - $ function $ \overline{\beta}_r = \overline{\beta}_r(m):  $

$$
 \overline{\beta}_r(m) := K_R(m) \cdot V_r(m) \cdot | \lambda_m(\cdot) |_{s, T, \nu},
$$
it is meaning $  \delta \in T, $
and suppose the finiteness of these function  for at last one value $   \overline{m} >  \overline{m}_0 := \max(2,p,q,s):  $

$$
\exists \overline{L}_r > \overline{m}_0 = \max(2,p,q,s) \ \Rightarrow  \overline{\beta}_r(\overline{L}_r) < \infty.  \eqno(6.13)
$$

\vspace{3mm}

\ We observe analogously to the proof of theorems 4.1a and 6.1 \\

\vspace{4mm}

{\bf Theorem 6.1a.} Let the condition (6.13) be satisfied. Then the (centered) r.p. $   \xi(t)   $ satisfies the CLT
in the space $ B^{o,p}_{\alpha,s} $ and moreover

$$
\sup_n \left| \left| \hspace{3mm} || S_n(\cdot) ||B^{o,p}_{\alpha,s} \hspace{3mm} \right| \right| G \overline{\beta}_r \le 1, \eqno(6.14)
$$
with correspondent tail estimation.\par

\vspace{4mm}

\section{ Non-asymptotical estimates.}

\vspace{3mm}

 \ We intend to obtain in this section the non-asymptotical estimates for the probability

$$
\overline{{\bf P}}(u) \stackrel{def}{=}
\sup_n {\bf P} \left( \ || S_n(\cdot) ||B^{o,p,q}_{\alpha,s} > u \ \right), \ u \ge e, \eqno(7.1)
$$
for example the exponential decreasing type bounds.\par

 \ Define a new function

$$
\kappa(m) := K_R(m) \cdot | d_m(t + z \delta, t)|_{p,T; \ q,T; \ s,T,\nu}, \ m > m_0 := \max(2,p,q,s) \eqno(7.2)
$$
and suppose its finiteness for the values $  m \in (m_0,  L'),  $ where $  L' \le \infty, $ i.e.  one can $ L' = \infty.$ \par
 We apply the inequality (5.8):

$$
\sup_n \left| \ || S_n(\cdot) ||B^{o,p,q}_{\alpha,s} \ \right|_{m,\Omega, {\bf P}} \le \kappa(m), \  m < L'.
 \eqno(7.3)
$$

\vspace{3mm}

{\bf Proposition 7.1.}\\

\vspace{3mm}

 Since  $ \tau \stackrel{def}{=} \ || S_n(\cdot) ||B^{o,p,q}_{\alpha,s} \ \in G\kappa $  and moreover $  ||\tau||G\kappa = 1, $
we conclude

$$
T_{\tau}(u) \le  \exp \left( - \tilde{\kappa}^*(\ln u) \right), \ u > e. \eqno(7.4)
$$
 where (we recall)  $ \tilde{\psi}(p) = p \cdot \ln \psi(p), $ and

$$
 \ g^*(y) = \sup_p (p |y| - g(p))
$$
is Young-Fenchel, or Legendre transform of the function $  g = g(p). $ \par

\vspace{3mm}

{\bf Example 7.1.}

\vspace{3mm}

 Suppose under the conditions of  remark 2.4. and hence under the conditions of  theorem 2.1

 $$
 d_m(t+ h,t) \le  C_m \ h^{\beta}, \ \beta  = \const \in (0,1].\eqno(7.5)
 $$

 \ If as before  $   C_m \le C \cdot m^{1/l}, \ l = \const > 0,   $

$$
\beta + 1/q > \alpha,  \hspace{6mm}  m \ge \max(2,p,q,s),
$$
then the r.p.  $ \xi(\cdot) $  satisfies the CLT in the space $ B^{o,p,q}_{\alpha, s}  $  and herewith

$$
 \sup_n T_{|| S_n ||B^{o,p,q}_{\alpha,s}}(u) \le  c_1(l; \alpha, \beta, p,q,s)  \ s^{-1/s} \  \times
$$

$$
 (\beta - \alpha + 1/q)^{-1/s} \cdot \exp \left\{-c_2(l;\alpha,\beta, p,q,s) \ u^{l/(l+1)} \right\}, \ u \ge e, \eqno(7.6)
$$
$  c_1, c_2 \in (0,\infty). $ \par

\vspace{3mm}

{\bf Example 7.2.}

\vspace{3mm}

 This time the variables $ m, \ C_m  $  can be considered as finite positive constants such that $ m \ge \max(2,p,q,s). $
Let again $ \beta + 1/q > \alpha.  $  We deduce on the basis of proposition  7.1

$$
 \sup_n T_{|| S_n ||B^{o,p,q}_{\alpha,s}}(u) \le  c_3(m; \alpha, \beta, p,q,s)  \ s^{-1/s} \  \times
$$

$$
 (\beta - \alpha + 1/q)^{-1/s} \ u^{-m}, \ u \ge 2. \eqno(7.7)
$$

\vspace{4mm}

{\bf Proposition 7.2.} \par

\vspace{3mm}

 \ The correspondent result for the classical Besov's spaces contains in theorem 6.1.  In detail:
 Let all the condition (6.8) be satisfied. Then the (centered) r.p. $   \xi(t)   $ satisfies the CLT
in the space $ B^{o,p}_{\alpha,s} $ and moreover

$$
\sup_n {\bf P} \left(|| S_n(\cdot) ||B^{o,p}_{\alpha,s} > u \right) \le
\exp \left( - \overline{\beta}^*(\ln u) \right), \ u > e. \ \eqno(7.8)
$$

\vspace{4mm}

\section{Concluding remarks. Multivariate case.}

\vspace{3mm}

 \hspace{6mm} {\bf A.} The possible generalizations may be undertaken in   two directions: multidimensional "time" $  t  $ and multivariate
increment  $ h, $ see detail definitions and investigations of these Besov's spaces
in \cite{Besov1}, \cite{Kufner1}, \cite{Milman1}, \cite{Simon1}. \par
 \ By our opinion, these generalizations are not hard for the  point of view of offered here problems. \par

\vspace{3mm}

 \ {\bf B.} The considered above Central Limit Theorem in Besov's spaces may be grounded also for (strong)
 stationary sequences $ \{  \xi_i  \},   $ superstrong mixingales, martingales etc., see \cite{Nachapetyan1}
\cite{Osekowski1}, \cite{Osekowski2}, \cite{Osekowski3}, \cite{Ostrovsky3}, \cite{Ostrovsky506}.\par
 \ In these articles was obtained in particular the analogs of Rosenthal's inequalities for
stationary sequences, superstrong mixingales and martingales. \par
 \ For example, the "Rosenthal's-Osekowski's" constant $ K_{R;Os}(p)  $ for the centered martingales, more exactly, for the
mean zero identical distributed martingale differences $  \{  \zeta_k \}, \ k = 1,2,\ldots  $  (relative arbitrary filtration),  i.e.

$$
K_{R; Os}(p) \stackrel{def}{=} \sup_{n \ge 1} \sup_{ \{\zeta_k\} } \left[ \frac{|n^{-1/2} \sum_{k=1}^n \zeta_k|_p}{|\zeta_1|_p} \right],
$$
may be estimated as follows

$$
K_{R; Os}(p)  \le 15.7858 \cdot \frac{p}{\ln p}, \hspace{5mm} p \ge 2,
$$
see \cite{Osekowski3}, \cite{Ostrovsky506}.\par

 \ For the  centered strictly stationary sequence of the r.v. $  \{  \zeta_k \}, \ k = 1,2,\ldots  $  satisfying the so-called superstrong
mixing condition with coefficient   $ \beta(k)  $ the analogous constant may be named as a constant of Rosenthal-Nachapetyan $  K_{R;N}(p)  $
\cite{Nachapetyan1} and may be estimated as follows

$$
 K_{R;N}(p)  \le C^p \cdot p^p \cdot \left[ \sum_{k=1}^{\infty} \beta(k) (k+1)^{ (p-2)/2 }  \right]^{1/p}.
$$

 \ See also \cite{Ostrovsky1},  chapter 2, section 2.9.\par


\vspace{4mm}

\begin{thebibliography}{99}

\vspace{4mm}

\bibitem{Adams1}
{\sc Adams R.A.} {\it  Anisotropic Sobolev Inequalities.}
Casopic pro Pestovani Matematiky, (Prague),  No. 3, 267—279.

\bibitem{Araujo1}
{\sc Araujo  A., Gine E. }
{\it The central limit theorem for real and Banach valued random variables.}
Wiley, (1980), London, New York.

\bibitem{Benedek1}
{\sc Benedek A. and Panzone  R.} {\it The space $  L_p $ with mixed norm.} Duke Math. J., {\bf 28}, (1961),  301-324.

\bibitem{Bennet1}
{\sc Bennett C. and Sharpley R.} {\it Interpolation of operators.}  Orlando, Academic Press Inc.,1988.

\bibitem{Besov1}
{\sc Besov O.V., Il’in V.P., Nikol’skii S.M. } {\it Integral representation of functions
and imbedding theorems.} Vol.1, 2; Scripta Series in Math., V.H.Winston
and Sons, (1979), New York, Toronto, Ontario, London.

\bibitem{Berg1}
{\sc J. Bergh and J. L\"ofstr\"om.} {\it Interpolation Spaces: an Introduction.} Springer, Berlin,
1976.

\bibitem{Billingsley1}
{\sc Billingsley P.} {\it Probability and measure.}
Wiley, 1979, London, New York.

\bibitem{Billingsley2}
{\sc Billingsley P.} {\it  Convergence of probability measures.}
Wiley, (1968), London, New York.


\bibitem{Boufoussi1}
{\sc B. Boufoussi, P. Chassaing and B. Roynette.} {\it A Kolmogorov criterion and an invariance
principle in Besov spaces.} Tech. Rept. No. 24, Prepublications de l'Institut
Elie Cartan, Nancy, 1993.

\bibitem{Buldygin1}
{\sc Buldygin, V.V. and Kozachenko, Yu.V.} {\it Metric characterization of random
variables and random processes.} Amer. Math. Soc., Providence, RI,  (2000).

\bibitem{Dudley1}
{\sc Dudley R.M.} {\it Uniform Central Limit Theorem}. Cambridge University Press, (1999)

\bibitem{Fageot1}
{\sc Julien Fageot,  Michael Unser,  John Paul Ward.}
{\it On the Besov Regularity of Periodic Lévy Noises.}\\
arXiv:1506.05740v1 [math.PR] 18 Jun 2015

\bibitem{Feneuil1}
{\sc Joseph Feneuil.} {\it Algebra properties for Besov spaces on unimodular Lie groups.}
arXiv:1505.06991v1 [math.AP] 26 May 2015

\bibitem{Fernique1}
 {\sc Fernique X.} (1975). {\it Regularite des trajectoires des
    function aleatiores gaussiennes. }  Ecole de Probablite de
    Saint-Flour, IV–1974, Lecture Notes in Mathematic. {\bf 480} 1–96, Springer Verlag, Berlin.

\bibitem{Fiorenza3}
 {\sc Fiorenza A., and Karadzhov G.E.} {\it Grand and small Lebesgue spaces and
       their analogs.} Consiglio Nationale Delle Ricerche, Instituto per le
      Applicazioni del Calcoto Mauro Picone, Sezione di Napoli, Rapporto tecnico n.
      272/03, (2005).

\bibitem{Frolov1}
{\sc Frolov A.S., Tchentzov N.N. } {\it On the calculation by the Monte-Carlo
method definite integrals depending on the parameters. } \\
Journal of Computational
Mathematics and Mathematical Physics, (1962), V. 2, Issue 4, p. 714-718 (in
Russian).

\bibitem{Gallagher1}
{\sc I. Gallagher and Y. Sire.} {\it Besov algebras on Lie groups of polynomial growth.}\\
 Studia Math., 212(2); 119- 139, 2012.


\bibitem{Garling1}
{\sc Garling D.J.H. }
{\it Functional Central Limit Theorems in Banach Spaces.}
 The Annals of Probability, Vol. 4, No. 4 (Aug., 1976), pp. 600-611

\bibitem{Gine1}
{\sc Gine E.} {\it On the Central Limit theorem for sample continuous processes.} Ann.
Probab. (1974), 2, 629-641.

\bibitem{Gine2}
{\sc Gine E., Zinn J.} {\it  Central Limit Theorem and Weak Laws of Large Numbers in certain Banach Spaces.  }
Z. Wahrscheinlichkeitstheory verw. Gebiete. {\bf 62}, (1983), 323-354.

\bibitem{Gogatishvili1}
{\sc Gogatishvili Amiran, Koskela Pekka,  Shanmugalingam Nageswari.} {\it Interpolation properties of Besov
spaces defined on metric spaces.} (English summary.) \\
 Math. Nachr. 283 (2010), no. 2, 215-231.

\bibitem{Grigorjeva1}
{\sc Grigorjeva M.L., Ostrovsky E.I.} {\it Calculation of Integrals on discontinuous
Functions by means of depending trials method.} Journal of Computational
Mathematics and Mathematical Physics, (1996), V. 36, Issue 12, p. 28-39 (in
Russian).

\bibitem{Hall1}
{\sc Hall P., Heyde C.C. } {\it Martingale Limit Theory and Applications.}  Academic
Press, New York. (1980)

\bibitem{Heinkel1}
{\sc Heinkel B.} Measures majorantes et le theoreme de la limite centrale dans
$ C(S). $ Z. Wahrscheinlichkeitstheory. verw. Geb., (1977). 38, 339-351.

\bibitem{Ibragimov1}
{\it Ibragimov R., Sharakhmetov Sh.}  {\it On an Exact Constant for the Rosenthal
Inequality.} \par
 Theory Probab. Applic., 42, N0 \ 2, 294-302, (1998).

\bibitem{Iwaniec2}
 {\sc Iwaniec T., P. Koskela P., and Onninen J.} {\it Mapping of finite distortion:
   Monotonicity and Continuity.}  Invent. Math. 144 (2001), 507-531.

\bibitem{Jain1}
{\sc Jain N.C. and Marcus M.B.} {\it Central limit theorem for $C(S)$ valued random
variables.} J. of Funct. Anal., (1975), 19, 216-231.

\bibitem{Johnson1}
{\sc Johnson W.B., Schechtman G., Zinn J.} {\it Best Constants in Moments Inequalities
for linear Combinations of independent and Changeable random Variables.} \\
 Ann. Probab., 1985, V. 13 pp. 234-253.

\bibitem{Kozachenko1}
 {\sc Kozachenko Yu. V., Ostrovsky E.I.} (1985). {\it The Banach Spaces of
      random Variables of subgaussian type.} Theory of Probab. and Math.
      Stat. (in Russian). Kiev, KSU, {\bf 32}, 43-57.

\bibitem{Yakovenko1}
{\sc Kozachenko, Yu.V. and Yakovenko, T.O.} {\it Conditions under which stochastic
processes belong to some function Orlicz spaces.} Bulletin of Kiyv University,
5, Kiev, (2002), 64-74.

\bibitem{Kufner1}
{\sc A. Kufner, O. John and S. Fu\"cik.} {\it Function Spaces. }  Noordhoff International Publishing, 1977.

\bibitem{Ledoux1}
 {\sc Ledoux M., Talagrand M.} (1991) {\it Probability in Banach Spaces.}\\
      Springer, Berlin, MR 1102015.


\bibitem{Liflyand1}
{\sc Liflyand E., Ostrovsky E., Sirota L.} {\it Structural Properties of Bilateral Grand Lebesgue Spaces.}
Turk. J. Math.; {\bf 34} (2010), 207-219.

\bibitem{Leoni1}
{\sc Leoni G. } {\it A first Course in Sobolev Spaces.} Graduate Studies in Mathematics,
v. 105, AMS, Providence, Rhode Island, (2009).

\bibitem{Lieb1}
{\sc Lieb E., Loss M.} {\it Analysis.} Providence, Rhode Island, 1997.


\bibitem{Marcos1}
{\sc Marcos M.A.} {\it Bessel potentials in Ahlfors regular metric spaces.}\\
arXiv:1506.08182v1 [math.CA] 26 Jun 2015

\bibitem{Milman1}
{\sc M.Milman and J.Xiao.} {\it The $  \infty \ - $  Besov capacity problem.}\\
arXiv:1506.01901v1 [math.AP] 5 Jun 2015

\bibitem{Morel1}
{\sc Morel Bruno.} {\it Weak convergence of summation processes in Besov spaces.}\\
Studia Matematica, 165, {\bf (1),} \ (2004), 19-38.

\bibitem{Nachapetyan1}
{\sc Nachapetyan B.S.} {\it On the certain criterion of weak dependence. } Probab. Theory Appl., (1980),
{\bf 2, } V. 26, 374-381.

\bibitem{Osekowski1}
{\sc A. Osekowski.} {\it Inequalities for dominated martingales.} Bernoulli 13 (2007), 54-79.

\bibitem{Osekowski2}
{\sc A. Osekowski.}  {\it Sharp martingale and semimartingale inequalities.} Monografie Matematyczne 72,
Birkh\"auser, 2012.

\bibitem{Osekowski3}
{\sc A. Osekowski.}  {\it  A Note on Burkholder-Rosenthal Inequality.} \\
 Bull. Polish Academy of Science, Math., 60, (2012), 177-185.

\bibitem{Ostrovsky1}
{\sc  Ostrovsky E.I.} (1999). {\it Exponential estimations for random Fields and its
Applications, (in Russian).}  Moscow-Obninsk, OINPE.

\bibitem{Ostrovsky2}
{\sc Ostrovsky E.I.} (1980).{\it On the support of probabilistic measures in separable
Banach spaces.} Soviet Mathematic, Doklady, v.255, No 6 pp. 836-838, (in Russian).

\bibitem{Ostrovsky301}
{\sc Ostrovsky E., L.Sirota L.}
{\it CLT for continuous random processes under approximations terms.}
arXiv:1304.0250v1 [math.PR] 31 Mar 2013

\bibitem{Ostrovsky3}
{\sc  Ostrovsky E. and Sirota L.} {\it  Moment and tail estimates for martingales and martingale transform,
with application to the martingale limit theorem in Banach spaces.}
arXiv:1206.4964v1 [math.PR] 21 Jun 2012

\bibitem{Ostrovsky302}
{\sc Ostrovsky E., Sirota L.} {\it Monte-Carlo method for multiple parametric integrals
calculation and solving of linear integral Fredholm equations of a second
kind, with confidence regions in uniform norm.}\\
 arXiv:1101.5381v1 [math.FA] 27 Jan 2011

\bibitem{Ostrovsky303}
{\sc  Ostrovsky E., Rogover E.} {\it Non - asymptotic exponential bounds for
MLE deviation under minimal conditions via classical and generic chaining methods.}\\
arXiv:0903.4062v1 [math.PR] 24 Mar 2009

\bibitem{Ostrovsky304}
{\sc Ostrovsky E., Sirota L.} {\it  Central Limit Theorem and exponnential tail estimations  in mixed (anisotropic)
Lebesgue spaces. }\\
arXiv:1308.5606v1 [math.PR] 26 Aug 2013

\bibitem{Ostrovsky305}
{\sc Ostrovsky E., Sirota L.} {\it Moment Banach spaces: Theory and Applications.}\\
HAIT Journal of Science and Engineering C, Volume 4, Issues 1-2, pp. 233-262.

\bibitem{Ostrovsky406}
{\sc Ostrovsky E., Sirota L.} {\it Factorable continuity of 4random fields,
with quantitative estimation.} \\
arXiv:1505.02839v1 [math.PR] 12 May 2015.


\bibitem{Ostrovsky407}
{\sc Ostrovsky E., Sirota L.} {\it  Central Limit Theorem and exponential Tail Estimations in mixed (anisotropic)
 Lebesgue Spaces. } \\
 arXiv:1308.5606v1 [math.PR] 26 Aug 2013

\bibitem{Ostrovsky408}
{\sc Ostrovsky E., Sirota L.} {\it Each Random Variable in separable Banach Space belongs to the Domain of Definition of
some inverse to compact linear non-random operator.}\\
arXiv:1410.3040v1 [math.PR] 12 Oct 2014

\bibitem{Ostrovsky301}
{\sc Ostrovsky E., L.Sirota L.}
{\it CLT for continuous random processes under approximations terms.}
arXiv:1304.0250v1 [math.PR] 31 Mar 2013

\bibitem{Ostrovsky603}
{\sc Ostrovsky E.} {\it Support of Borelian measures in separable Banach spaces. }
arXiv:0808.3248v1 [math.FA] 24 Aug 2008

\bibitem{Ostrovsky502}
{\sc  Ostrovsky E. and Sirota L.}  {\it  Schl\"omilch and Bell series for Bessel's functions, with
probabilistic applications.}
arXiv:0804.0089v1 [math.CV] 1 Apr 2008

\bibitem{Ostrovsky506}
{\sc  Ostrovsky E. and Sirota L.} {\it Sharp moment estimates for polynomial martingales. }
arXiv:1410.0739v1 [math.PR] 3 Oct 2014

\bibitem{Peetre1}
{\sc J. Peetre.} {\it New Thoughts on Besov Spaces.} Duke Univ. Math. Ser. 1, Duke Univ.,
Durham, 1976, 58-77.

\bibitem{Pinelis1}
{\sc Pinelis I.F. and Utev S.A.} {\it Estimates of Moment of sums of independent
random variables.} \\
 Theory Probab. Appl., 1984, V. 29 pp. 574-578.

\bibitem{Pisier1}
{\sc  Pisier G,, J. Zinn.} {\it On the limit theorems for random variables with values in the spaces } $ L_p, \ 2  \le p  < \infty. $
 Z. Wahrscheinlichkeitstheorie verw. Gebiete 41, 289-304, (1978).

\bibitem{Pisier2}
{\sc Pizier G.}  {\it Condition d’entropic assupant la continuite de certain processus et
applications a l’analyse harmonique.} Seminaire d’analyse fonctionnalle. (1980),
Exp. 13,  p. 23-24.

\bibitem{Prokhorov1}
{\sc Prokhorov Yu.V.} {\it Convergense of Random Processes and Limit Theorems
of Probability Theory.} Probab. Theory Appl., (1956), V. 1, 177-238.

\bibitem{Ratckauskas1}
{\sc Rackauskas A, Suquet Ch.}
{\it Central limit theorems in H\"ölder topologies for Banach space valued random fields.}
Teor. Veroyatnost. i Primenen., 2004, Volume 49, Issue 1, Pages 109–125 (Mi tvp238)

\bibitem{Ratchkauskas2}
{\sc  Ratchkauskas A., Ch. Suquet.} {\it Necessary and sufficient condition for the H\"olderian functional central limit theorem. }
J. Theoret. Probab. 17 (2004) 221–243.

\bibitem{Ratchkauskas3}
{\sc Ratchkauskas A, Suquet Ch.} {\it H\"older norm test statistics for epidemic change. } J. Statist. Plann. Inference, {\bf 126}, (2004), 495-520.

\bibitem{Ratchkauskas4}
{\sc Ratchkauskas A, Suquet Ch.} {\it Testing epidemic changes of infinite dimensional parameters.}\\
 Stat. Inference Stoch. Process,  {\bf 9},  (2006), 111-134.

\bibitem{Rosenthal1}
{\sc Rosenthal H.P.} {\it On the subspaces of $ L_p \ (p > 2) $ spanned by sequences of
independent Variables. } \\
 \ Israel J. Math., 1970, V.3 pp. 273-253.

\bibitem{Simon1}
{\sc J. Simon.} {\it Sobolev, Besov and Nikolskii fractional spaces: imbeddings and comparisons for vector valued
spaces on an interval.} \\
 \ Annali di Matematica pura ed applicata (IV), Vol. LCVII (1990), 117-148.

\bibitem{Talagrand1}
 {\sc Talagrand M.} (1996). {\it Majorizing measure: The generic chaining.}
 Ann. Probab., {\bf 24} 1049-1103. MR1825156

\bibitem{Talagrand2}
 {\sc Talagrand M.} (2005). {\it The Generic Chaining. Upper and
     Lower Bounds of Stochastic Processes.} Springer, Berlin. MR2133757.

\bibitem{Triebel1}
{\sc  Triebel H.} {\it Theory of Function Spaces.}  Birkh\"auser, Basel, 1992.

\bibitem{Yakovenko2}
{\sc Yakovenko Tetyana.} {\it  Stochastic processes in some Besov spaces. }\\
Theory of Stochastic Processes,
Vol.13, {\bf 29}, no.1-2, 2007, pp.308-315.

\bibitem{Zinn1}
{\sc Zinn J.  } {\it A Note on the Central Limit Theorem in Banach Spaces.}
 Ann. Probab. Volume 5, Number 2 (1977), 283-286.


\end{thebibliography}
\end{document}